\documentclass[11pt]{amsart}
\usepackage {amsmath, amssymb, a4wide, epsfig, enumerate, psfrag}
\usepackage[latin1]{inputenc}
\date{\today}

\keywords{}
 \author{Romain Dujardin}

 \title[Continuity of Lyapunov exponents]{Continuity of Lyapunov exponents for 
 polynomial automorphisms of $\cd$}

\subjclass[2000]{Primary: 37F45, Secondary: 37F10, 32H50}

\address{Universit{\'e} Paris 7 et Institut de Math{\'e}matiques de Jussieu,
         {\'E}quipe G{\'e}om{\'e}trie et Dynamique,
         Case 7012, 2  place Jussieu,
         75251 Paris Cedex 05,
         France.}
\email{dujardin@math.jussieu.fr}

\newcommand{\cc}{\mathbb{C}}
\newcommand{\bb}{\mathbb{B}}

\newcommand{\dd}{\mathbb{D}}

\newcommand{\nn}{\mathbb{N}}
\newcommand{\e}{\varepsilon}
\newcommand{\cv}{\rightarrow}

\newcommand{\fr}{\partial}
\newcommand{\om}{\Omega}
\newcommand{\set}[1]{\left\{#1\right\}}
\newcommand{\norm}[1]{\left\Vert#1\right\Vert}
\newcommand{\abs}[1]{\left\vert#1\right\vert}
\newcommand{\cd}{\cc^2}
\newcommand{\pd}{{\mathbb{P}^2}}
\newcommand{\pu}{{\mathbb{P}^1}}
\newcommand{\rest}[1]{ \arrowvert_{#1}}

\newcommand{\unsur}[1]{\frac{1}{#1}}
\newcommand{\el}{\mathcal{L}}
\newcommand{\qq}{\mathcal{Q}}
\newcommand{\rond}{\hspace{-.1em}\circ\hspace{-.1em}}
\newcommand{\ms}{\mathrm{sm}}

\DeclareMathOperator{\supp}{Supp}

\DeclareMathOperator{\mult}{mult}

\newtheorem{prop} {Proposition} [section]
\newtheorem{thm}[prop] {Theorem} 
\newtheorem{defi}[prop] {Definition}
\newtheorem{lem}[prop] {Lemma}
\newtheorem{cor}[prop]{Corollary}

\newtheorem{defprop}[prop]{Definition-Proposition}
\theoremstyle{remark}

\newtheorem{rmk}[prop]{Remark}

\begin{document}

\begin{abstract}
We prove two continuity theorems for the Lyapunov exponents of the
maximal entropy measure of polynomial automorphisms of $\cd$. 
The first continuity result holds for any family
of polynomial automorphisms of constant dynamical degree. 
The second result is the continuity of the upper exponent for families
degenerating to a 1-dimensional map. 
\end{abstract}

\maketitle

\section{Introduction}\label{sec:intro}

\subsection{The main results}
Let $(f_\lambda)_{\lambda\in \Lambda}$ 
be a holomorphic family of polynomial diffeomorphisms of $\cd$, parameterized by a
complex manifold $\Lambda$, and with 
 constant dynamical  degree $d= \lim_{n\cv\infty}
(\deg(f^n))^{1/n}>1$. 
Then, according to the work of Bedford,
Lyubich and Smillie \cite{bs3, bls}, for each $\lambda$, $f_\lambda$
admits a natural invariant measure $\mu_\lambda$, with two non zero
exponents of opposite signs $\chi^-(f_\lambda)<0<\chi^+(f_\lambda)$, 
which is the unique measure of maximal entropy $\log d$.  Our first
main result is the following.

\begin{thm} \label{thm:main}
Let $(f_\lambda)_{\lambda\in\Lambda}$ 
be a holomorphic family of  polynomial
 automorphisms of $\cd$ of dynamical degree $d>1$, parameterized by a complex
  manifold $\Lambda$. 
Then the Lyapunov exponents of the maximal entropy measure are
 continuous functions of $\lambda$.
\end{thm}

Bedford and Smillie \cite{bs3} had previously shown that $\lambda\mapsto 
\chi^+(f_\lambda)$ is
plurisubharmonic (psh for short), hence in particular upper
semicontinuous.  Since polynomial automorphisms of $\cd$ have constant Jacobian,
the sum 
\begin{equation}\label{eq:sum}
\chi^+(f_\lambda)+ \chi^-(f_\lambda) =
 \log\abs{\mathrm{Jac}(f_\lambda)}
\end{equation} is
a pluriharmonic function of $\lambda$, so $\chi^-$ is
plurisuperharmonic.

As a consequence of Young's formula \cite{young} for the dimension of
$\mu$, 
\begin{equation}
\label{eq:young}
\dim(\mu)= \log d \left(\unsur{\chi^+}-\unsur{\chi^-}\right),
\end{equation}
we obtain the following corollary.

\begin{cor}
The Hausdorff dimension of $\mu_\lambda$ is a continuous function on
parameter space.
\end{cor}

The reader may object that the statement of Theorem \ref{thm:main}
would sound more natural if the family $(f_\lambda)$ was only supposed
to depend continuously on $\lambda$. Actually this is not more
general:  the space of polynomial diffeomorphisms of degree $d$
is a (reducible) finite dimensional complex variety, so there are
finitely many holomorphic families whose union  covers {\em all} 
polynomial diffeomorphisms of fixed degree.

\medskip

An important issue in the theory of polynomial automorphisms of the
(real or complex) plane is the study of their degenerations to one
dimensional maps. The sample model is the family of  {\em
 generalized  H{\'e}non mappings}
$$g_{p,b}:(z,w)\longmapsto (p(z)-bw, z),$$  degenerating to the
1-dimensional map
$z\mapsto (z,p(z))$ when the Jacobian $b$ tends to zero. Here $p$ is a
polynomial in $z$.
Usually, a dynamical assumption (like hyperbolicity, Collet-Eckmann, or
renormalizability) is made on $p$, and it is studied how
this assumption influences the dynamics of the 2 dimensional mapping 
$g_{p,b}$ when $b$ is small. 

It is convenient for our purposes to conjugate $g_{p,b}$ (by a linear
map) so that it 
rewrites as $(z,w)\mapsto (aw+p(z), az)$ and the degenerate form
becomes $(p(z),0)$. We work in the following general setting (see \S
\ref{sec:extended} for some examples): 
consider a complex manifold $\Lambda$, and a family of polynomial mappings of $\cd$,
depending holomorphically on $\lambda$,  of the 
form 
\begin{equation}\label{eq:degen2}
f_\lambda(z,w)= (p(z),0) + R_\lambda(z,w), 
\end{equation}
with $p$ monic and of degree $d$, and $R_\lambda$ a polynomial mapping
of $\cd$, of degree $\leq d-1$, vanishing identically for
$\lambda=\lambda_0$. We assume that for $\lambda$ outside some
hypersurface $\Lambda_{\rm degen}$, $f_\lambda$ is an automorphism of
$\cd$.

For $\lambda\in \Lambda_{\rm
  degen}$ close to zero, $f_\lambda$ reduces to a
 1 dimensional mapping of degree $d$, so it has a unique measure of maximal
entropy $\log d$ (see Lemma \ref{lem:gather2}), 
and we still denote by $\chi^+(f_\lambda)$  its unique (positive)
Lyapunov exponent.

 Our continuity result reads as follows. Notice that we make no
 hypothesis on $p$.

\begin{thm}\label{thm:bis}
Assume that $(f_\lambda)$ is a degenerating family of polynomial
automorphims of the form \eqref{eq:degen2}. Then $\lambda\mapsto 
\chi^+(f_\lambda)$ is  continuous at $\lambda_0$.  
\end{thm} 

As before, upper semicontinuity has been established in
\cite{bs3}, so we only have to prove lower semicontinuity. 

Observe
 further that by \eqref{eq:sum}, $\chi^-(f_\lambda)$ tends to $-\infty$ as
$\lambda\cv \lambda_0$. Hence by combining
  Young's formula \eqref{eq:young} and the formula of Ma{\~n}e \cite{mane}
$\dim(\mu_p)=\log d / \chi^+(p)$, we conclude that 
under the assumptions of the previous theorem, $\lambda \mapsto
\dim(\mu_\lambda)$
is  continuous at $\lambda_0$. 

\subsection{One dimensional maps}
To understand our motivations and methods, it is useful to make a d{\'e}tour through
one dimensional dynamics. Let $(f_\lambda)_{\lambda\in\Lambda}$ be any 
holomorphic family of rational maps of $\pu$, of degree $d$. For every
$\lambda$, $f_\lambda$ admits a unique measure of maximal entropy
$\log d$ \cite{lyub, flm} with a unique positive exponent
  $\chi^+(f_\lambda)\geq \log d/2$. This defines a natural psh function
  $\lambda\mapsto \chi^+(f_\lambda)$ on parameter space. It is a
  theorem by Ma{\~n}e  \cite{mane}  that this function is continuous. 

The importance of this function is underlined by the following
remarkable connection between analysis and dynamics: the bifurcation
locus of $(f_\lambda)$ is precisely the support of the current $dd^c
(\chi^+(f_\lambda))$ \cite{demarco2}. The geometry of this current and
its successive powers is explored in \cite{bas-ber, df}, where the
continuity of the potential plays an important technical role.  

\medskip

There are basically two ways of proving continuity in this
setting.  In a
``direct" approach, Ma{\~n}e \cite{mane} considers 
the following integral formula for $\chi^+$:
\begin{equation}\label{eq:naive}
\chi^+(f_\lambda)= \int_\pu \log\norm{Df_\lambda} d\mu_\lambda,
\end{equation}
where $\norm{\cdot}$ is the operator norm associated to some Riemannian
metric on $\pu$. It is an easy fact that $\lambda\mapsto \mu_\lambda$
is continuous in the weak topology. This does not immediately imply
the continuity of $\chi^+$ because of the possible critical points of
$f_\lambda$ on $\supp(f_\lambda)$. 
To overcome this difficulty, it is necessary to study the uniform
integrability of functions with logarithmic singularities, which
follows for instance by standard potential theoretic estimates. This
method has been generalized for the sum of exponents of polynomial-like
mappings in higher dimension by Pham \cite{pham}. 

The second approach uses more elaborate formulas for $\chi^+$. The
prototypical example is the following well known  formula
(due to Manning \cite{manning}, Przytycki \cite{prz} and Sibony
\cite{sib-ucla}) for 
the Lyapunov exponent  when $f_\lambda$ is a polynomial 
\begin{equation}\label{eq:prz}
\chi^+(f_\lambda) = \log d + \sum_{c \text{ critical}} G_\lambda
  (c),
\end{equation} where $G_\lambda$ is the dynamical Green function (see \S
  \ref{sec:extended} for more details). 
The H{\"o}lder continuity of $\chi^+$ follows immediately. 
This has been generalized to  rational maps of $\pu$ by DeMarco
  \cite{demarco2} and to certain sums of exponents in higher dimension
  in \cite{bj, bas-ber}. 

\subsection{Methods} Let us return to polynomial automorphisms of
  $\cd$. There is an integral formula analogous to \eqref{eq:naive}
  for $\chi^+(f_\lambda)$:
$$\chi^+(f_\lambda) = \int \log\norm{Df_\lambda\rest{e^u} } d\mu_\lambda, $$
where $e^u$ is the {\em unstable direction}, provided by the
Osedelets Theorem, which is defined almost everywhere
 and varies measurably (see \cite{bls} for an adapted
presentation). It is unclear how to study the
variation of this data as a function of $\lambda$.

On the other hand there is a formula, due to Bedford and Smillie
\cite{bs5},  analogous to \eqref{eq:prz} in this context. 
We do not state this formula here (see \S \ref{subs:bs5} for
details), but we indicate that it involves
{\em unstable critical points}, that is,  critical points of the Green
function $G^+$ restricted to unstable manifolds. 
This seems to be an important notion regarding the geometry of the
unstable lamination, but it is still poorly understood.  

\medskip

The main idea of the proof of Theorem \ref{thm:main} is to prove that
the formula of \cite{bs5} varies lower semicontinously. For this, we
use a result of \cite{connex} on the geometry of the
 unstable lamination for unstably disconnected mappings (see
 \S\ref{subs:structure} for the definition of unstable
 disconnectedness) 
that allows a precise counting of unstable critical points.

Likewise, in the proof of Theorem \ref{thm:bis}, we study the
``convergence'' of the formula of \cite{bs5} to \eqref{eq:prz} when the
Jacobian tends to zero. This involves a description of the
geometry of the unstable lamination for perturbations of 1-dimensional
maps, regardless of hyperbolicity. 

We see that, besides the intrinsic interest of Theorems \ref{thm:main}
and \ref{thm:bis}, another interesting point in the paper is the
contribution to the understanding of the unstable lamination of
polynomial diffeomorphisms. 

\subsection{Outline} The structure of the paper is as
follows: sections \ref{sec:main} and \ref{sec:extended} are
respectively devoted to the proofs of Theorems \ref{thm:main} and
\ref{thm:bis}. 

As opposite to the 1-dimensional situation, one feature of unstable
critical points is that they are invariant under the dynamics.  In an
Appendix, we briefly develop a notion of {\em fastest rate of escape} for 
  critical points, generalizing the corresponding 1-dimensional notion. This is a
natural level (of the Green function) where looking at unstable
critical points, that provides a new dynamically defined function on
parameter space. We also show that it gives an upper estimate for the
Lyapunov exponent (Theorem \ref{thm:appendix}).


\section{Proof of theorem \ref{thm:main}}\label{sec:main}

\subsection{Basics}\label{subs:basics}
We will use some standard facts
 from the dynamics of polynomial automorphisms of $\cd$. General
 references are \cite{fm, bs1, bls, sib}. As usual $K^+$
 (resp. $K^-$) is the set of points with bounded forward
 (resp. backward) orbits. We denote the invariant currents  by
 $T^\pm$ and the Green functions by $G^\pm$ so that $T^\pm=dd^cG^\pm$.. 

We also need the notion of {\em uniformly laminar current}. Recall that
a positive (closed) current is uniformly laminar  if locally, up to change of
coordinates, it is  a current of
integration over a measured family of graphs over the unit disk.
The reader is referred to \cite{bls, structure} for more details, and
to \cite{de} for basics on positive currents.

\medskip 

We now introduce a  notion of {\em horizontality} of objects 
relative to a projection. 

\begin{defi}
Let $\om$ be an open set, together with a
holomorphic locally trivial fibration
 $\varphi: \om\cv D$, where $D$ is a disk in
$\cc$.   We say that an 
analytic subset or current is {\em horizontal} relative to
$\varphi$ if  $\supp(T)$ 
is (locally uniformly) relatively  compact in the fibers of
$\varphi$. More specifically, we ask that for every $\zeta\in D$, there
exists a neighborhood $D'\ni \zeta$ and a compact subset $K$ of $\om$
such that for every $\zeta'\in D'$, $\varphi^{-1}(\zeta')\cap \supp(T)
\subset K$. 
\end{defi}

The simplest case is when 
$\om=D\times \cc$ or $D\times D$, and $\varphi$ is the
first projection. In this case we just say that $T$ is horizontal in
$D\times \cc$ (resp. $D\times D$). Some basic properties of horizontal
currents in this context can be found in \cite{sl, hl} (see also \cite{ds}). 

It is clear that if $V$ is a horizontal analytic subset, the number
(couting multiplicities) of
intersection points of $V$ and the fibers of $\varphi$ is locally
(hence globally) constant. We call this number the {\em
  degree} of $V$ relative to $\varphi$. When there is no danger of
confusion, we drop the mention to $\varphi$ and denote it by 
$\deg(V)$. 

There is a related statement for currents. For simplicity we
only consider the case where $T$ admits a global potential in $\om$. 

\begin{defprop}
If $T=dd^c u$ is a horizontal positive closed current
  relative to some $\varphi:\om\cv D$, then for every $\zeta\in D$, 
 the wedge product 
$T\wedge [\varphi^{-1}(\zeta)]$ is well defined, and  its
  total mass does not depend on $\zeta$. By definition, we call this mass
  the {\em slice mass} of $T$, $\ms(T)$. 
\end{defprop}

 For instance, when $V$ is a
horizontal analytic set, $\ms([V])=\mathrm{deg}(V)$. 

\begin{proof} Pick  any fiber $F=\varphi^{-1}(\zeta)$. The potential 
 $u$ is pluriharmonic near $\fr F$, so $u$ cannot be
 identically $-\infty$ on any connected component of $F$. Hence
 $u\rest{F}\in L^1_{\rm loc}(F)$ and $dd^cu\wedge [F]$ is well
 defined. 

Let $F_0\Subset F$ be an open subset with smooth (oriented) 
boundary $\gamma$ 
such that $\supp(T)\cap F\subset F_0$.  Then by Stokes' Theorem $\int_F dd^c u =
\int_\gamma d^cu$. By assumption, the form $d^cu$ is closed near
$\gamma$, so by Stokes' Theorem again,  the integral 
$\int_\gamma d^cu$ is insensitive to small perturbations of $\gamma$ in
other fibers. We conclude that $\int_{\varphi^{-1}(\zeta)} dd^c u$ is
  locally constant as a function of $\zeta$.
\end{proof}

\begin{rmk}\label{rmk:appendix}
To get the conclusion if the proposition, 
we do not really need $\varphi$ to be a fibration. What is  important is the 
existence of a smoothly varying $\gamma$ enclosing $\supp(T)$ in the fibers, so that 
we can apply Stokes' Theorem in the same way. This observation
will be useful in the proof of Theorem  \ref{thm:appendix}. 
\end{rmk}


\subsection{Preparation}\label{subs:prepar}
Let $(f_\lambda)_{\lambda\in\Lambda}$ be a holomorphic family of
polynomial automorphisms of $\cd$  of
constant  dynamical degree $d= \lim_{n\cv\infty} (\deg(f^n))^{1/n} >1$. 
By the  work of Friedland and Milnor \cite{fm},
a polynomial automorphism of $\cd$ has non trivial dynamics if and
only if its dynamical degree is larger than 1. 
Without loss of generality, we may
assume that $\Lambda$ is an open subset of $\mathbb{C}^N$ for some
$N$, and $0\in \Lambda$. We denote the upper
Lyapunov exponent by  $\chi^+(f_\lambda)$  ($\chi^+(\lambda)$ for short). 
We  will prove the continuity of
the upper Lyapunov exponent $\chi^+$
 at $\lambda=0$. 

In the next lemma we show that in the neighborhood of $0$, the
mappings $f_\lambda$ may be written in a simpler form. We classically
consider the extension of polynomial automorphisms as rational maps of
the projective plane $\pd$, and denote by $I^+(f)$ (resp $I^-(f)$)
the indeterminacy point of $f$ (resp $f^{-1}$); see \cite{sib} for more details. 
We also consider homogeneous coordinates $[z:w:t]$ on $\pd$ so that our
$\cd$ imbeds as $(z,w) \mapsto [z:w:1]$.  

\begin{lem}\label{lem:prepar}
There exists an open subset $\Lambda'\subset\Lambda$ containing 0, and
a holomorphic family $(\varphi_\lambda)_{\lambda\in \Lambda'}$
of polynomial automorphisms of $\cd$ such that  for $\lambda \in
\Lambda'$, 
$$\varphi_\lambda^{-1} f_\lambda \varphi_\lambda(z,w)= (z^d, 0)+
\text{ lower degree terms}.$$   
 \end{lem}
\begin{proof} It is well known \cite{fm} that every polynomial automorphism 
of $\cd$ with non trivial dynamics  is conjugated in the group of
  polynomial automorphisms of $\cd$ to a composition 
generalized complex H{\'e}non  mappings, 
which can be written under the form  $(z,w)\mapsto
  (aw+p(z), az)$, with $p$  monic and of degree $d$. 
Let $\varphi_0$ be such a conjugating map for $f_0$. We get that $I^+(
\varphi_0^{-1}f_0\varphi_0)=[0:1:0]$ and $I^-(
\varphi_0^{-1}f_0\varphi_0)=[1:0:0]$. 

Now for small $\lambda$, $I^+(\varphi_0^{-1}f_\lambda\varphi_0)$ is close to
$[0:1:0]$. Also, $f^{-1}_\lambda$ depends holomorphically on
$\lambda$. Indeed, we know that for each $\lambda$, $f^{-1}_\lambda$
is a polynomial map of degree $d$, and it is of course completely
determined by  its restriction to  a small open set. In such an open
set, the holomorphic dependence of $f^{-1}_\lambda$ on $\lambda$
follows from the Implicit Function Theorem. In particular we get that 
$\lambda\mapsto I^-(\varphi_0^{-1}f_\lambda\varphi_0)$ is holomorphic. 

Hence for small $\lambda$, $I^+(\varphi_0^{-1}f_\lambda\varphi_0)
\neq I^-(\varphi_0^{-1}f_\lambda\varphi_0)$ and we get that 
$\varphi_0^{-1}f_\lambda\varphi_0$  has non trivial
dynamics. Conjugating with a holomorphically varying linear map
$\ell_\lambda$, with $\ell_0= id$ 
we can ensure that
$I^+((\varphi_0 \ell_\lambda)^{-1}f_\lambda(\varphi_0
\ell_\lambda))= [0:1:0]$ and $I^-((\varphi_0 \ell_\lambda)^{-1}f_\lambda(\varphi_0
\ell_\lambda))= [1:0:0]$. 

Inspecting  the higher order terms  shows that 
$$(\varphi_0 \ell_\lambda)^{-1} f_\lambda(\varphi_0
\ell_\lambda)(z,w)= (c(\lambda)z^d, 0) + \text{ lower degree terms},$$
with $c(0)=1$. For small $\lambda$, there exists a holomorphic
determination of $c(\lambda)^{-1/(d-1)}$, and conjugating 
with $(z,w)\mapsto (c(\lambda)^{-1/(d-1)}z,w)$ gives the desired
result.
\end{proof}

Lyapunov exponents are preserved under smooth
conjugacy, so, slightly abusing notation, in the following we replace 
$f_\lambda$ by $\varphi_\lambda^{-1} f_\lambda\varphi_\lambda$ and
$\Lambda$ by $\Lambda'$, so that we can assume that the conclusion of
the lemma holds. 

\medskip

Since polynomial automorphisms have constant Jacobian, 
$\chi^+(f)+\chi^-(f) = \log\abs{\mathrm{Jac}(f)}$
varies continuously. Hence  it is enough to consider the continuity of
$\chi^+$.
Bedford and Smillie  proved in \cite{bs3} that $\lambda\mapsto
\chi^+(\lambda)$ is plurisubharmonic, hence
in particular upper semicontinuous. 
We  recall their  argument for completeness.  
Let $\norm{\cdot}$ be any operator norm on the space of linear maps of
$\cd$. For $n\in \nn$, let
$$\chi^+_n(\lambda)=\unsur{n} \int \log \norm{Df_\lambda^n}d\mu_\lambda,$$ 
where $Df$ is the tangent map. From the chain rule we get that
$\chi^+_n$ satisfies the following subadditivity property
$$(m+n)\chi^+_{m+n}\leq m\chi^+_m+ n\chi^+_n.$$ In particular the
sequence $\chi^+_{2^n}$ is decreasing. 

On the other hand, for every $n$, $\lambda \mapsto\chi^+_n(\lambda)$ is
continuous. This is a property of continuous variation of the maximal
entropy measure, which itself easily follows from the joint continuity
of the escape rate functions $G^\pm_\lambda(z,w)$ in the variable 
$(\lambda,(z,w))$ \cite{bs1}. We conclude that $\chi^+$ is a
decreasing limit of continuous functions, hence upper semicontinuous. 

Once upper semicontinuity is established, plurisubharmonicity follows
from instance from the fact that $\chi^+$ can be evaluated on saddle
orbits \cite{bls2}.

\medskip

We infer that it is enough to prove that $\chi^+$ is lower
semicontinuous, that is,
\begin{equation}\label{eq:sci}
\liminf_{\lambda\cv 0}\chi^+(\lambda)\geq
\chi^+(0).
\end{equation}
Another result from \cite{bs3} is that $\chi^+$ is 
always bounded from below by $\log d$, so \eqref{eq:sci} 
is automatic when $\chi^+(0)=\log d$.
 Moreover, replacing $f_0$ by $f_0^{-1}$ if necessary, we may
assume that $\abs{\mathrm{Jac}(f_0)}\leq 1$, in which case we have the following
 \cite{bs5,bs6} 
$$\chi^+(f_0)=\log(d)\Leftrightarrow J(f_0) \text{ is connected }
\Leftrightarrow f_0 \text{ is unstably connected.}$$
We will explain what unstable connectedness means shortly.
For the moment we 
conclude that Theorem \ref{thm:main} reduces to the following proposition.

\begin{prop}\label{prop:main}
Assume that $(f_\lambda)_{\lambda\in \Lambda}$ is a holomorphic family 
of polynomial automorphisms of $\cd$ of degree $d$, 
with $I^+(f_\lambda)=[0:1:0]$ and $I^-(f_\lambda)=[1:0:0]$. 
Assume further that $\abs{\mathrm{Jac}(f_0)}\leq 1$ and $f_0$ is unstably 
disconnected. 
Then $\chi^+$ is lower semicontinuous at $f_0$.
\end{prop}


\subsection{Laminar structure}\label{subs:structure}
In this paragraph we explain some results of \cite{connex} on the laminar structure of 
$T^-$ for unstably disconnected mappings.

Let $f$ be a polynomial automorphism satisfying the conclusion of Lemma 
\ref{lem:prepar}. Then if $R$ is large enough, $f$ stretches the horizontal direction 
in $D_R^2$ (here $D_R=D(0,R)$) in the sense that in $D_R^2$
the image of a horizontal manifold is horizontal --the degree is multiplied by $d$. 

The following definition was introduced in  \cite{bs6}.

\begin{defi}\label{def:unstdisc}
$f$ is {\em unstably disconnected} if for some saddle point $p$, $W^u(p)\cap K^+$
has a compact component (for the topology induced by the isomorphism
$W^u(p)\simeq \cc$).
\end{defi}

This condition is actually independent of the saddle point $p$ \cite{bs6} (see also 
\cite[Prop. 1.8]{connex}). Now if $V_0$ is a disk in $W^u(p)$ such that 
$V_0\cap K^+\Subset V_0$, pushing $V_0$ sufficiently
 many times gives rise to a horizontal 
disk $V$ of finite degree in $D_R^2$, lying inside $W^u(p)$. 

In \cite[Prop. 2.3]{connex} 
we proved that $V$ is {\em  subordinate to} $T^-$ in the sense that there 
exist a non trivial uniformly laminar current $S\leq T^-$, made up of a family of 
disjoint horizontal disks of degree $\mathrm{deg}(V)$, extending $V$. 
As an easy consequence, 
we get the following basic result on the structure 
of $T^-$ \cite[Th. 2.4]{connex}.

\begin{thm}\label{thm:structure}
If $f$ is unstably disconnected, then  
there exists a sequence of  currents $T^-_k$ in $D_R^2$, with 
 $$T^-\rest{D_R^2}=\sum_{k=1}^\infty T_k^-,$$
and  $T^-_k$ is the integration current over a family of disjoint 
horizontal disks of degree $k$.

\end{thm}

Of course the same holds for any sub-bidisk of the form $D\times D_R$,
but we stress that 
{\em the decomposition depends on the bidisk}. 
Notice that, even if the degree is bounded,  this need not imply that  $T^-$ 
is  uniformly
laminar. Indeed a sequence of disks of degree 2
may accumulate on a disk of degree 1, giving rise to some folding.

\begin{defi}
We say that a horizontal current $T$ in a bidisk
has finite degree $K$ if it is an integral of disjoint 
horizontal analytic subsets of degree $0\leq k\leq K$. These analytic subsets 
are said to be subordinate to $T$. 
\end{defi}

We denote by $\mathcal{F}_K(D\times \cc)$ the
class of currents of finite degree $K$. This is relative to some
projection which should be clear from the context.

We will need the following (expected) proposition. This does not seems
to be a direct 
consequence of the Choquet Representation Theorem because of the
disjointness assumption. 

\begin{prop}\label{prop:choquet}
Assume that $(T_n)$ is a sequence of currents in
$\mathcal{F}_k(D\times\cc)$, weakly converging to $T$, and whose
supports are contained in a fixed vertically compact subset. Then
$T\in \mathcal{F}_k$. 
\end{prop}

\begin{proof} For the proof, we say that two 
subvarieties are {\em compatible} if
they do not have isolated intersection points. 
We work in the unit bidisk $\bb=\dd^2$.
Assume that $T_n$ is a sequence of (uniformly) horizontal
currents of finite degree $K$, converging to $T$. Using Bishop's compactness theorem
for curves of bounded volume  \cite{bishop}
and the Hurwitz lemma, we will first construct a family of 
horizontal analytic sets to which $T$ should be subordinate, and 
then indeed prove the subordination by using the Hahn-Banach Theorem. 

\medskip

The first observation is that if $V_n$ is a sequence of curves in
$\bb$, with $V_n\subset \dd\times D_ {1-\e}$ for some $\e$, and of
degree $\leq K$ , then it has locally finite area 
(see e.g.  \cite[Lemma 3.6]{sl}), so
there exists a converging subsequence by Bishop's Theorem. It is
obvious that any cluster value must have degree $\leq K$. Notice that
the limit may be singular. Another useful observation is that if
$(V_n)$ and $(V_{n}')$ are two converging sequences of varieties of
degree $\leq K$, with $V_n\cap V_n'=\emptyset$, then by Hurwitz'
Theorem, their limits are compatible. 

\medskip

Let $X_0=\limsup (\supp (T_n))\supset \supp(T)$, and $x_0\in X_0$. There exists a
sequence $x_n\in \supp(T_n)$ converging to $x_0$; let $V_n(x_n)$ be
the analytic set  subordinate to $T_n$ passing through $x_n$. There exists a
subsequence $\varphi_0(n)$ so that
$V_{\varphi_0(n)}(x_{\varphi_0(n)})$ 
converges to a $V(x_0)\ni x_0$ of degree $\leq K$.

Let now $X_1= \limsup (\supp (T_{\varphi_0(n)}))\supset \supp(T)$. By
construction, $x_0\in X_1$. Let $r_1$ be the supremum of the radii of
balls centered at points of $X_1$ and avoiding $x_0$, and let $x_1\in
X_1$ so that $x_0\notin B(x_1,\frac{r_1}{2})$. Let $x_{\varphi_0(n)}\in 
\supp (T_{\varphi_0(n)})$, $x_{\varphi_0(n)}\cv x_1$, and 
$V_{\varphi_0(n)}(x_{\varphi_0(n)})$ be subordinate to
$T_{\varphi_0(n)}$ through $x_{\varphi_0(n)}$. 
 From $\varphi_0(n)$, 
extract a subsequence $\varphi_1(n)$ such that this sequence of
analytic sets converges to $V_1\ni x_1$. By construction, $V_0$ and $V_1$
are compatible, and have degree $\leq K$. 

Inductively we construct a sequence of successive extractions
$\varphi_k(n)$, a decreasing sequence of closed sets $X_k \supset
\supp(T)$, a sequence of points $x_k$ with $x_0, \ldots, x_k \in X_k$,
together with a family of compatible analytic subsets $V_0, \ldots,
V_k$.  Let then $X_{k+1} = \limsup (\supp (T_{\varphi_k(n)}))$.  The
point $x_{k+1}$ is chosen in the following way: let $r_{k+1}$ be the
supremum of radii of balls centered at points of $X_{k+1}$ and
avoiding $x_0,\ldots, x_k$. We choose $x_{k+1}$ such that $x_0,
\ldots, x_k \notin B(x_{k+1},\frac{r_{k+1}}{2})$, and the extraction 
$\varphi_{k+1}$, and $V_{k+1}$ as for the case $k=0$.

Let $X= \bigcap_{k\geq 0} X_k \supset \supp(T)$. It is an
exercise to show that the sequence $(x_k)$ is dense in
$X$. Attached to each $x_k$ there is an analytic set $V_k$ of degree
$\leq K$. Let $\el$ be the closure of the family $(V_k)$; $\el$  is a family
of compatible analytic sets  of degree $\leq K$. In particular for every
$x\in X$, there is a unique $V(x)\in \el$ containing $x$.   

\medskip

It remains to see that $T$ is an integral of the varieties in  $\el$,
or equivalently, that $T\in \mathrm{Conv}(\el)$. Assume not. Then by the
Hahn-Banach Theorem, there exists a smooth test (1,1) form $\phi$ such
that $\langle
T,\phi\rangle < 0$ while $\langle [V],\phi\rangle > 0$ for every
$V\in\el$. Let $\varphi_\infty$ be a diagonal extraction of the
$\varphi_k$, and for ease of notation denote by $\varphi_\infty(n)$ by
$n'$. 
For large $n$', $\langle T_{n'},\phi\rangle < 0$, so there
exists an analytic set $V_{n'}$ subordinate to $T_{n'}$ such that 
$\langle [V_{n'}],\phi\rangle < 0$. Extract a further subsequence 
(still denoted by $n'$) so
that $V_{n'}$ converges to some $V$ ($V$ needn't be contained in
$\supp(T)$, this is the reason for the limsup above). We claim that
$V\in \el$. Indeed, let $x_{n'}\in V_{n'}$ be  convergent to $x$. Then
$x\in X$. Now if $V\neq V(x)$, since near $x$, $V\cap V(x)=\set{x}$, 
 $V$ would have isolated intersection
points with $V_{n'}$ for large $n'$, a contradiction. 

By construction, $V\in \el$ and $\langle [V],\phi\rangle \leq 0$ which
is contradictory. This finishes the proof.
\end{proof}


\subsection{The formula of \cite{bs5}}\label{subs:bs5}
In this paragraph we describe the formula of \cite{bs5} for the
Lyapounov exponents. The results in this section do not require
unstable disconnectedness.  We start with a temporary definition.

\begin{defi}\label{def:temporary}
An unstable critical point is a critical point of $G^+\rest{W^u(p)}$,
where $p$ is some saddle periodic point.
\end{defi}

Every unstable manifold is equidistributed along the unstable current
$T^-$, so, following \cite{bs5}, it is natural to define an unstable
critical measure by extending the notion of critical point to Pesin
unstable manifolds and integrate against the transverse measure. The
formal definition is a bit delicate because of the weakness of the laminar
structure of $T^-$. Here we provide a precise definition using the results of 
\cite{structure}.

An embedded holomorphic  disk $\Delta$ is said to be {\em subordinate} to $T^-$ if
there exists a uniformly laminar current $0<S\leq T^-$ such that
$\Delta$ lies inside a leaf of the lamination induced by $S$. 
We can now revise Definition \ref{def:temporary}.

\begin{defi}
An unstable critical point is a critical point of $G^+\rest{\Delta}$,
where $\Delta$ is any disk subordinate to $T^-$.
\end{defi}

By  {\em flow box} we mean a piece of lamination, biholomorphic to a union of
graphs in the bidisk. In \cite[Th. 1.1]{structure} we proved that if $\el$ is
any flow box, $T^-$ induces by restriction an invariant transverse
measure on $\el$. Moreover it is clear from the construction of the
laminar structure of $T^-$ \cite{bls,bs5} that there exists a 
countable family $(\el_i)$ of (overlapping) flow boxes so that every disk
subordinate to $T^-$ is contained in some union of $\el_i$. We say
that $(\el_i)$ is a {\em complete system} of flow boxes associated to
$T^-$. 

Observe that by the Radon Nikodym Theorem, 
it is possible to define the supremum of two measures in the following
way: if $\mu_1$ and $\mu_2$ are two $\sigma$-finite positive measures,
there exist measurable bounded functions $f_i$, $i=1,2$ such that 
$\mu_i=f_i(\mu_1+\mu_2)$. By definition  $\sup(\mu_1, \mu_2) =
\sup(f_1,f_2)(\mu_1+\mu_2)$. We may hence define the supremum of a
finite family of measures by induction, and  if 
$(\mu_i)_{\set{i\geq 1}}$ is  sequence of positive
measures, we define $\sup(\mu_i)_{\set{i\geq 1}}$ as the increasing limit of 
 $\sup(\mu_i)_{\set{1\leq i\leq I}}$ as $I\cv\infty$. 
 It is not clear of course that the limiting measure  will have locally finite
mass. 

We are now ready to define the critical measure.

\begin{defi}\label{def:critmeas}
Let $\el$ be a flow box $\el = \bigcup_{t\in \tau}{L_t}$, where $\tau$
is a global transversal to $\el$,
 and write $T^-\rest{\el}= \int_\tau
 [L_t]d\mu_\el(t)$, where $\mu_\el$ is the measure induced by $T^-$ on
 $\tau$. 

The critical measure restricted to  $\el$ is defined by 
$$\mu^-_c\rest{\el} = \int_\tau
\left[\mathrm{Crit}(G^+\rest{\el_t})\right] d\mu_\el(t),$$
where $\mathrm{Crit}(G^+\rest{\el_t})$ is the sum of point masses at
critical points of $G^+$ on $\el_t$, counting multiplicities. 

The global critical measure $\mu_c^-$ is now defined by $\mu_c^- =
\sup\left(\mu_c\rest{\el_i}\right)$ where $(\el_i)$ 
is a countable  complete system of flow boxes associated to
$T^-$, as above.
\end{defi} 

This definition is independent of the choice of
the system of flow boxes. Indeed if $(\el'_j)_j$ is another choice, the
respective leaves of $(\el_i)_i$ and $(\el'_j)_j$ are compatible because
$T^-\wedge T^-=0$ (see also \cite[Th. 1.1]{structure}), and the result
follows by
considering the  system of flow boxes $(\el_i\cap \el'_j)_{i,j}$. 
 
\medskip

We can now state the main result of \cite{bs5}, upon which the proof
of continuity will be based. 

\begin{thm}[\cite{bs5}]\label{thm:bs5}
The upper Lyapounov exponent of the maximal entropy measure satisfies 
\begin{equation}\label{eq:bs5}
\chi^+ = \log d + \int_{1\leq G^+<d} G^+d\mu_c^- = 
 \log d + \int_P G^+d\mu_c^-,
\end{equation} where $P$ is any measurable
 fundamental domain for $f\rest{\cd\setminus K^+}$.
\end{thm}


\subsection{The total mass of the critical measure}
In this paragraph we consider an unstably disconnected 
 polynomial automorphism $f$, satisfying the conclusion of Lemma
 \ref{lem:prepar}.  We will give a formula for the mass of the
 critical measure in certain domains. 

We introduce the B{\"o}ttcher function $\varphi^+$.  If $R$ is a positive
real number, we classically
denote by $V_R^+$ the forward invariant open set
$$V_R^+=\set{(z,w), \ \abs{z}> R, \abs{w}<\abs{z}}.$$ 
Recall that $f(z,w)= (z^d, 0) +\ l.o.t.$ near
infinity. By analogy with the B{\"o}ttcher coordinate in one variable
dynamics, for  large enough $R$ and $(z,w)\in V_R^+$, we can define \cite{ho}
$$\varphi^+(z,w)= \lim_{n\cv\infty } \left(\pi_1\rond
  f^n\right)^{\unsur{d^n}}$$
  where the $(1/d^n)^{\rm th}$  root is chosen so that
  $\varphi^+(z,w) = z+ O(1)$ at infinity. It is a holomorphic function
  in $V_R^+$ satisfying the functional equation $\varphi^+\rond f =
  (\varphi^+)^d$ and $G^+$ equals  $\log\abs{\varphi^+}$.

Geometrically, it should be understood as an invariant 
first projection near infinity. 
The following lemma is obvious.

\begin{lem}
Unstable critical
  points in $V_R^+$ are points of tangency between unstable manifolds
  and the fibers of $\varphi^+$. Moreover the multiplicity of a
  critical point  as a
  critical point of $G^+$  and as a vertical tangency coincide. 
\end{lem}

Since we are not going to consider $\varphi^-$, from now on we write
$\varphi$ for $\varphi^+$. For the same reason 
we drop the minus sign from $\mu_c^-$. 
By Rouch{\'e}'s Theorem, $(\varphi, w)$ is a coordinate system in $V_R^+$
for large enough $R$. Furthermore, if $Q$ is a bounded simply connected 
 open subset of  $\set{\abs{z}>R+C}$, where $C$ is such that
 $\abs{\varphi(z,w)-z}<C$ in $V_R^+$, then 
 $T^-$ is horizontal in 
  $\varphi^{-1}(Q)\cap V_R^+$, and admits a
 decomposition like the one in Theorem \ref{thm:structure}, relative to
 $\varphi$. Indeed, consider
 such a decomposition in $\set{\abs{z}, \abs{w}<R'}$, with $R'\gg R$
 and  restrict it to $\varphi^{-1}(Q)\cap V_R^+$. 

Another important remark is that since the fibers of $\varphi$ are
vertical graphs in $V_R^+$, $T^-$ has slice mass $1$ with respect to the
projection $\varphi$.
We now fix such a $R$, and for notational simplicity, we write
$\varphi^{-1}(Q)$ for 
$\varphi^{-1}(Q)\cap V_R^+$. 

The decomposition of Theorem \ref{thm:structure} induces  a formula
for the mass of $\mu_c$ in $\varphi^{-1}(Q)$. 

\begin{prop}\label{prop:mass}
Let $Q$ be a bounded simply connected 
 open subset of  $\set{\abs{z}>R}$, and $T^-\rest{\varphi^{-1}(Q)}
 = \sum_{k=1}^\infty T_k$ be the decomposition of $T^-$ in 
$\varphi^{-1}(Q)$ relative to $\varphi$, as  above. Then 
$$ \mu_c(\varphi^{-1}(Q)) = 
\sum_{k=1}^\infty \frac{k-1}{k} \ms(T_k).$$
\end{prop}

Since $\sum \ms(T_k)=1$
this implies in particular (without using Theorem \ref{thm:bs5}) that the critical
measure has locally finite mass. 

\begin{proof}
It is enough to compute the contribution of $T_k$ for each $k$. Recall
that $T_k$ is made of horizontal disks of degree $k$ over $Q$. Consider
such a disk $\Delta$ 
and a uniformly laminar current $S\leq T_k$ supported in a small
tubular neighborhood of $\Delta$. Fix a vertical fiber
$\varphi^{-1}(z)$ transverse to $\Delta$, hence intersecting
$\Delta$ in exactly $k$ points. If $\supp(S)$ is close enough to
$\Delta$, the same holds for every leaf of $S$. Hence, as a flow box,
the transverse measure of $S$ is $\unsur{k}\ms(S)$. Now, by the
Riemann-Hurwitz formula, each leaf of $S$ carries exactly $k-1$
critical points, counting multiplicities, so the contribution of $S$
to the mass of the critical measure is
$\frac{k-1}{k}\ms(S)$. 

Exhausting $T_k^-$ by such uniformly laminar currents finishes the
proof.
\end{proof}

\subsection{Conclusion}\label{subs:proof}

We return to the setting of \S \ref{subs:prepar}, and consider a
family of polynomial automorphisms  satisfying the assumptions
of Proposition \ref{prop:main}. Under these assumptions, it is an easy
fact that the escape radius $R$ of the previous paragraph 
is locally uniformly bounded, and the resulting B{\"o}ttcher function
$\varphi_\lambda$ is defined on a fixed $V_R^+$ and 
depends holomorphically on $\lambda$. 

We claim that Proposition \ref{prop:main} follows from the following.

\begin{prop}\label{prop:sci} Under the above assumptions, 
let $Q$ be a simply connected open subset of
$\set{\abs{z}>R}$, with piecewise smooth boundary. 
Then the critical mass of $\varphi_\lambda^{-1}(Q)$
is lower semicontinuous, that is  (with obvious notation)
$$\liminf_{\lambda\cv 0} \mu_{c, \lambda}\left(\varphi_\lambda^{-1}(Q)
\right) \geq \mu_{c, 0}\left(\varphi_0^{-1}(Q)
\right).$$
\end{prop}

\begin{proof}[Proof of Proposition \ref{prop:main} assuming
  Proposition \ref{prop:sci}]
By \eqref{eq:bs5}, it is enough to prove that $$\lambda\longmapsto\int_{\set{A\leq
    G^+_\lambda <dA}}G^+_\lambda d\mu_{c,\lambda}$$ is lower semicontinuous at
$\lambda=0$ for some $A$. We choose $A$ so large that
if $Q\subset \set{e^A<\abs{z}<e^{dA}}$, then $T^-$ is horizontal in
$\varphi^{-1}(Q)\cap V_R^+$. Slightly moving $A$ if necessary we can
further assume that 
 $\int_{\set{G_0^+=A}}G^+_0 d\mu_{c,0}=0$. 

Consider a sequence $(\mathcal{Q}_n)$ of subdivisions of the annulus 
$\set{e^A<\abs{z}<e^{dA}}$ by  simply
connected and piecewise smoothly bounded pieces (``squares") of size
smaller than
$\unsur{n}$.  We may further assume that
$\mu_{c,0}\left(\varphi_0^{-1}(\fr \qq_n)\right)=0$. Then, 
since $G_0^+$ is continuous and constant along the fibers of 
$\varphi_0$, 
$$  \int_{\set{A\leq
    G^+_0 <dA}}G^+_0 d\mu_{c,0} = \lim_{n\cv\infty}
\sum_{Q\in\qq_n} \left( \inf_{\varphi_0^{-1}(Q)} G_0^+\right) \
\mu_{c,0}\left(\varphi_0^{-1}(Q)\right),$$ where the limit on the
    right hand side is increasing. For $\lambda\neq 0$, we have a
    similar result, except that the critical measure could charge the
    boundary of some subdivision, so that  
$$  \int_{\set{A\leq
    G^+_\lambda <dA}}G^+_\lambda d\mu_{c,\lambda} \geq  \lim_{n\cv\infty}
\sum_{Q\in\qq_n} \left(\inf_{\varphi_\lambda^{-1}(Q)} G_\lambda^+\right) \
\mu_{c,\lambda}\left(\varphi_0^{-1}(Q)\right).$$
Let $h_n(\lambda)$ be the sum of the right hand side, and
    $h_\infty(\lambda)$ be its (increasing) limit. Notice further that 
since $G^+_\lambda= \log\abs{\varphi_\lambda}$,
$\inf_{\varphi_\lambda^{-1}(Q)} (G_\lambda^+) = \inf_{z\in Q}
    (\log{\abs{z}})$ does not depend on $\lambda$. Hence by Proposition
    \ref{prop:sci}, $h_n(\lambda)$ is lower semicontinuous at 0 for
    every $n$. Lower semicontinuity is preserved under increasing
    limits so we get that 
$$\liminf_{\lambda\cv 0} \int_{\set{A\leq
    G^+_\lambda <dA}}G^+_\lambda d\mu_{c,\lambda} \geq 
\liminf_{\lambda\cv 0} h_\infty(\lambda) \geq  h_\infty(0) = 
\int_{\set{A\leq
    G^+_0 <dA}}G^+_0 d\mu_{c,0}, $$ and the result follows.
\end{proof}

\begin{proof}[Proof of Proposition \ref{prop:sci}]
We fix $Q$ as in the statement of the proposition, and denote by
$Q^\delta= \set{z\in Q, \ \mathrm{dist}(z,\fr Q)>\delta}$. For small
$\delta$, $Q^\delta$ is a topological disk.

The following lemma is easy and left to the reader.

\begin{lem} \label{lem:horiz}
If $\delta>0$ is fixed, then there exists a neighborhood $N$ of $0\in
\Lambda$ depending only on $\delta$ 
such that if  $\lambda\in N$ and $C$ is  any horizontal curve  of degree $k$ in
$\varphi^{-1}_\lambda(Q)$ (relative to the projection $\varphi_\lambda$), 
then 
$C\cap \varphi^{-1}_0(Q^\delta)$ is a horizontal curve of
degree $k$ relative to $\varphi_0$.
\end{lem}

Observe that if $C$ is a disk, $C\cap \varphi^{-1}_0(Q^\delta)$ can be
disconnected, so by the maximum principle it is a union of disks. 

Recall the notation $\mathcal{F}_K(\varphi_\lambda^{-1}(Q))$ for the
set of currents of finite degree $K$ over $Q$ relative to the
projection $\varphi_\lambda$. 
By the previous lemma, if $\delta$ is fixed and $\lambda$ is small
enough, then   $$S\in \mathcal{F}_K(\varphi_\lambda^{-1}(Q))
\Rightarrow S\rest{\varphi^{-1}_0(Q^\delta)}\in
\mathcal{F}_K(\varphi_0^{-1}(Q^\delta)).$$ An important further remark
is that the slice mass of $S$ is invariant under small perturbations
of the transversal,
so in particular it does not depend on $\lambda$. 

\medskip

We will prove that for every $\delta>0$, 
\begin{equation}\label{eq:delta}
 \liminf_{\lambda\cv 0} \mu_{c, \lambda}\left(\varphi_\lambda^{-1}(Q)
\right) \geq \mu_{c, 0}\left(\varphi_0^{-1}(Q^\delta)
\right),
\end{equation}
whence the desired result by letting $\delta\cv0$. For this, we use
decompositions $T^-(\lambda)= \sum T_k(\lambda)$ with the following
conventions:
\begin{itemize}
\item[-] for $\lambda\neq 0$, the decomposition is relative to
  $\varphi^{-1}_\lambda(Q)$ (and the projection $\varphi_\lambda$);
\item[-] for $\lambda=0$, 
the decomposition is relative to
  $\varphi^{-1}_0(Q^\delta)$ (and the projection $\varphi_0$).
\end{itemize}
We also use the following notation:
 $m_k(\star)= \ms(T_k(\star))$, $T^-_{\leq K}(\star)=
\sum_{1\leq k \leq K} T_k(\star)$, and $M_{\leq K}(\star)
=\ms (T_{\leq K}(\star))$, where $\star$ stands for 0 or $\lambda$. 

The first observation is that the locus of unstable disconnectivity is
open (see \cite[\S 2.1]{connex}, this  actually follows from the
discussion below). In particular for small $\lambda$, Proposition
\ref{prop:mass} applies. 

As $\lambda\cv 0$, by Proposition \ref{prop:choquet}
every cluster value $S$ of $T^-_{\leq K}(\lambda)$ in
$\varphi_0^{-1}(Q^\delta)$ belongs to
$\mathcal{F}_K(\varphi_0^{-1}(Q^\delta))$ and satisfies $S\leq
T^-(0)$, so $S\leq T^-_{\leq K}(0)$. 
From this we get that 
\begin{equation}\label{eq:1}
\limsup_{\lambda\cv 0} M_{\leq K}(\lambda) \leq M_{\leq K}(0).
\end{equation}
We want to prove that 
\begin{equation}\label{eq:2}
\liminf_{\lambda\cv 0} \sum_{k=1}^\infty \frac{k-1}{k} m_k(\lambda)
\geq \sum_{k=1}^\infty \frac{k-1}{k} m_k(0).
\end{equation}
Since $T^-(\star)$ has total slice mass 1, we infer that 
$\sum_{k=1}^\infty  m_k(\star)=1$. In particular, 
 \begin{equation}\label{eq:3}
\sum_{k=1}^\infty \frac{k-1}{k} m_k(\star) = 1- \sum_{k=1}^\infty 
\frac{m_k(\star)}{k} .
\end{equation}
We now make the following classical ``integration by parts" on
series:
$m_k(\star) = M_{\leq k}(\star)-M_{\leq k-1}(\star)$, with the
convention that $M_0=0$, so that 
$$\sum_{k=1}^\infty  \frac{m_k(\star)}{k}  = \sum_{k=1}^\infty
\frac{M_{\leq k}(\star)}{k(k+1)}.$$ So by \eqref{eq:1} we get that 
$$\limsup_{\lambda\cv 0}
\sum_{k=1}^\infty  \frac{m_k(\lambda)}{k}  \leq 
\sum_{k=1}^\infty  \frac{m_k(0)}{k} ,$$
which from \eqref{eq:3} implies \eqref{eq:2}. This concludes the proof.
\end{proof}


\section{Extended parameter space}\label{sec:extended}

In this section we prove that the upper Lyapunov exponent  is
continuous when the $(f_\lambda)$ degenerate to a 1-dimensional
map. In a preliminary subsection, we start with some general
 considerations on degenerating families of polynomial automorphisms.

\subsection{Degenerating families of polynomial automorphisms}
There is no  classification of degenerating families 
  of polynomial automorphisms in the literature, so  we
  will consider a general, but possibly non exhaustive, situation. 

Consider as  parameter space a neighborhood $\Lambda$ of the
origin in $\cc^n$, and a family of polynomial mappings of $\cd$,
depending holomorphically on $\lambda$,  of the 
form 
\begin{equation}\label{eq:degen}
f_\lambda(z,w)= (p(z),0) + R_\lambda(z,w),
\end{equation}
with  $p$ monic and of degree $d$, $R_\lambda$ a polynomial mapping of degree 
$\leq d-1$ with $R_0\equiv 0$,
and such that for $\lambda$ outside some hypersurface $\Lambda_{\rm
  degen}$, $f_\lambda$ is an automorphism of $\cd$.
  
A typical  such situation is as follows: 
 fix integers $d_1, \ldots, d_m$, with $d_i\geq 2$, and let 
$d=d_1\cdots d_m$. Consider
a family of polynomial automorphisms
\begin{equation}\label{eq:henon}
f_\lambda = f_{m, \lambda}\rond\cdots\rond f_{1, \lambda}\text{ ,
  where } 
f_{i,\lambda} = (a_{i,\lambda}w+ p_{i,\lambda}(z),a_{i,\lambda}z),
\end{equation}
where
 $a_{i,\lambda}$ are holomorphic functions of $\lambda$,  
$p_{i,\lambda}$ are monic polynomials of degree $d_i$, depending
holomorphically on $\lambda$, and $a_{i_0,0}=0$ for some $i_0$. By
  conjugating with an appropriate 
composition of the $f_i$,  without loss of generality 
we may assume that $i_0=m$. After this
  conjugacy, this family of automorphisms is of the form
  \eqref{eq:degen}. 

Recall from \cite{fm} that  any polynomial
  automorphism of $\cd$ with non trivial dynamics is of the form
  \eqref{eq:henon}, up to conjugacy. 

\medskip

In the next lemmas, we collect some well known results on the 
extension of the dynamics to $\Lambda_{\rm degen}$. Recall that 
$V^+_R= \set{(z,w), \ \abs{z}>\abs{w}, \ \abs{z}>R}$,
and  $\pi_1$ denotes the first projection in $\cd$.

\begin{lem}\label{lem:gather}
For $\lambda\in \Lambda$, let $G^+_\lambda(z,w)= \lim_{n\cv\infty}
\log^+\norm{f_\lambda^n(z,w)}$. Then $G^+$ is continous on $\Lambda$ 
as a function of $(\lambda, z,w)$. 

Also, reducing $\Lambda$ if necessary, there exists a fixed $R>0$ so
that the B{\"o}ttcher function $\varphi_\lambda(z,w) = \lim_{n\cv\infty} 
(\pi_1\rond f^n(z,w))^\unsur{d^n}$ is well defined in $V_R^+$, and is
jointly holomorphic in $(\lambda, z,w)$.
\end{lem}

\begin{proof} We sketch the proof for completeness. Extending the
  polynomial mappings $f_\lambda$ to the projective plane $\pd$ (with
  coordinates $[z:w:t]$), we find that the indeterminacy set
  $I_\lambda^+$ is  constant and equal to $[0:1:0]$. Also, $f_\lambda$
  maps the line at infinity with $I^+$ deleted on $[1:0:0]$ so when
  $\lambda\notin \Lambda_{\mathrm{degen}}$, the indeterminacy set of
  $f_\lambda^{-1}$ is $[1:0:0]$. In particular $f_\lambda$ is regular
  in the sense of Sibony \cite{sib}, and has entropy $\log d$. 

There are several ways of proving the continuity of $(\lambda,
z,w)\mapsto G^+_\lambda(z,w)$. Following  \cite{sib}, we use the
homogeneous lift $F_\lambda$ of $f_\lambda$ to $\cc^3$ (see \cite{sib}
for more details). Let $\pi:\cc^3\cv\pd$ be the natural
projection. For $\lambda$ close to 0, there exists uniform constants 
$c$ and $C$ such that if $\pi(p)\in \pd$ is far away from
$I^+$, $c\norm{p}^d\leq F_\lambda(p) \leq C\norm{p}^d$. Moreover,
every point escapes
any small neighborhood of $I^+$ after finitely many iterations of
$f_\lambda$, whenever $\lambda$ is in $\Lambda_{\rm degen}$ or not. 
Indeed when $\lambda \notin \Lambda_{\rm degen}$ this is classical, 
and when $\lambda\in \Lambda_{\rm degen}$, this follows from the fact
that the image of $f_\lambda$ is a variety not going through $I^+$. 
 We
conclude that the limit defining $G^+_{\lambda}(z,w)$ is locally
uniform. 

The second assertion of the lemma is a simple consequence of the
construction of the B{\"o}ttcher function
 \cite{ho} and is left to the reader.
\end{proof}

A consequence which will be important for us is 
that when $\lambda=0$, both $G^+_0$ and $\varphi_0$ depend only
on $z$, so that at the limit the unstable  critical points become points where 
 unstable manifolds have vertical tangencies. 

\begin{lem}\label{lem:gather2}
If $\lambda\in \Lambda_{\rm degen}$ is close enough to zero, then
$f_\lambda(\cd)=M_\lambda$ is a (possibly singular) 
subvariety of degree $\leq d$, and $f_\lambda\rest{M_\lambda}$  is
conjugated to a polynomial map of degree $d$.  

Moreover $T^-_\lambda:= [M_\lambda]/\deg(M_\lambda)$ is the unique
closed positive current invariant under $\unsur{d}(f_\lambda)_*$, and 
the maximal entropy measure of $f_\lambda$ is $T_\lambda^+\wedge
T_\lambda^-= dd^cG_\lambda^+\wedge T_\lambda^-$. 
\end{lem}

Recall that by definition a subvariety is irreducible.
We do not know any example where $M_\lambda$ is singular.
As a simple 
illustration, let 
$$g_{a,b}(z,w)= (aw+q(z),az)\circ (bw+r(z),bz).$$
This is a degenerating family of polynomial automorphisms of the form
\eqref{eq:degen}. 
Here the parameter space is $\Lambda= \cd_{a,b}$,  $\Lambda_{\rm
  degen} = \set{ab=0}$, $p=q\rond r$ and
$d=\deg(q)\deg(r)$. If $b=0$ and $a\neq 0$, $M_\lambda = 
M_{a,0}=\set{(q(t),
  at), t\in \cc}$ has degree $\deg(q)$. 

\begin{proof}
 Assume $\lambda\in \Lambda_{\rm degen}$ and  is close to 
 zero. We first prove that $f_\lambda(\cd)=M_\lambda$ 
is a  subvariety of $\cd$ close to $\set{w=0}$. It is convenient to
 consider the meromorphic extension of $f_\lambda$ to $\pd$.

For $\lambda\in \Lambda_{\rm degen}$, $f_\lambda$ is not an
automorphism. Since $f_\lambda$ is approximated by automorphisms, it
has constant Jacobian. If the Jacobian was a non zero
constant, $f_\lambda$ would be  a local diffeomorphism of finite
degree. Now if for some $q\in \cd$, $f_\lambda^{-1}(q)$ was a
finite set with at least two elements, this would persistently hold in
the neighborhood of $\lambda$ in $\Lambda$, which is not possible. We
conclude that for $\lambda\in \Lambda_{\rm degen}$, $f_\lambda$ has
zero Jacobian. 

Consider a generic line $\overline L$ in $\pd$ not going through
the indeterminacy set of $f_\lambda$. Then $f_\lambda(\overline L)$ is 
an irreducible analytic subset of degree $\leq d$ in $\pd$. 
Let $L=\overline L \cap \cd$. Since a generic $L$ must
intersect every fiber of $f_\lambda$, except possibly finitely many of
them, we conclude that
$f_\lambda(\cd)=f_\lambda(L) = M_\lambda$ is irreducible and of degree
$\leq d$. 

Furthermore, $f_\lambda\rest{L}: L\cv M_\lambda$ gives rise to a
finite-to-one map $\cc\cv M_\lambda$. If $M_\lambda$ is smooth, this
directly imply that $M_\lambda$ is isomorphic to $\cc$. If not, we
only get that the desingularization of   $M_\lambda$ is isomorphic to
$\cc$.

\medskip

We compute the topological degree of $f_\lambda\rest{M_\lambda}$. Fix 
a large bidisk $V_R$ such that every point of $M_\lambda\setminus D_R^2$ escapes to 
infinity, and $M_\lambda$ is horizontal in $D_R^2$.  A generic 
vertical line $L^v$  in $D_R^2$ intersects 
$M_\lambda$ in $\mathrm{deg}(M_\lambda)$ points.  $f^{-1}_\lambda (L^v)\cap 
D_R^2$ is a vertical curve of degree $d$, thus intersecting $M_\lambda$ at 
$d\mathrm{deg}(M_\lambda)$ points, that are the preimages of the previous ones 
under $f_\lambda\rest{M_\lambda}$. It follows that the degree of
$f_\lambda\rest{M_\lambda}$ is $d$.

We conclude that $f_\lambda$ induces a
holomorphic self map of (the desingularization of) $M_{\lambda}$ 
 of degree $d$, that is, a polynomial of degree $d$ on $\cc$. 

\medskip

The fact that $[M_\lambda]/\deg(M_\lambda)$ is the only current invariant under 
$\unsur{d}f_*$ is obvious. And by the functional equation $G^+_\lambda
\rond f_\lambda = dG^+_\lambda$, we infer that
$dd^c(G^+_\lambda\rest{M_\lambda})$ is a non atomic measure of constant Jacobian
$d$ relative to  $f_\lambda$, hence the unique measure of maximal entropy.
\end{proof}
  
\subsection{Continuity of the upper exponent}  

Whenever degenerate or not, for small 
$\lambda$, $f_\lambda$ has a unique measure of
  maximal entropy $\log d$, with only one positive exponent, still denoted 
by $\chi^+(f_\lambda)$.

For convenience we recall the statement of the continuity result. 

\begin{thm}\label{thm:bisbis}
Assume that $(f_\lambda)$ is a degenerating family of polynomial
automorphims of the form \eqref{eq:degen}. Then $\chi^+(f_\lambda)$ converges to
$\chi^+(p)$ as $\lambda\cv 0$. 
\end{thm} 

\begin{proof} As before, upper semicontinuity has been established in
\cite{bs3}, so we only prove lower semicontinuity. 
Again, the tools will be
the formula \eqref{eq:bs5} of \cite{bs5}, and the Manning-Przytycki formula for
$\chi^+(p)$  \cite{manning, prz}
\begin{equation}\label{eq:manning}
\chi^+(p) = \log d  + \sum_{c \text{ critical}} G_p(c).
\end{equation} Of course,
only escaping critical points contribute to the sum.

The starting point  is to rewrite the Manning-Przytycki
 formula in the light of the formula of 
Bedford and Smillie. Let $p$ be  monic and of degree
$d$, and define the critical measure $\mu_{c,0}$ for $p$ as follows
(the subscript $0$ is used because $p=f_0$) 
$$\mu_{c,0}= \sum_{c \text{ critical escaping}} \mult_p(c) \sum_{k\geq 0} 
 \unsur{d^k} \delta_{f^k(c)},$$
where $\mult_p(c)$ is the multiplicity of $c$ as a critical point
    (e.g. $\mult_{z^2}(0)=1$). 

Let $G_{\rm max}$ be the maximum of $G_p(c)$ over all critical
points. By the invariance relation for $G_p$,
 we immediately get the following rewriting of \eqref{eq:manning}:
\begin{equation}\label{eq:manningbis}
\forall A\geq G_{\rm max}, \ \chi^+(p) = \log d + \int_{A\leq G_p <
  dA} G_p d\mu_{c,0}.
\end{equation}
We will have to
understand how unstable critical points of $f_\lambda$ 
degenerate to the escaping postcritical points of $f_0$. 

\medskip
 First, a brief outline of the proof, which will comprise
  several sublemmas. We start with the easier situation where $\lambda$
  approaches 0 along the hypersurface $\Lambda_{\rm degen}$ of
  degenerate parameters (Proposition 
  \ref{prop:degen}).  
In a second step, we give an interpretation
  of the critical measure for $p$ in terms of counting ramified and
  unramified inverse branches of some domain $Q$ outside $K_p$ (Lemma
  \ref{lem:branch} and Corollary \ref{cor:branch}). 
The third step is to make the
  connection with 2-dimensional dynamics, by inspecting the geometry
  of iterated submanifolds, as projected along  the B{\"o}ttcher
  fibration (Lemmas \ref{lem:tang} and \ref{lem:tub}). We then conclude
  the proof by using the decomposition of $T^-$ given by Theorem
  \ref{thm:structure}.
\medskip

\noindent {\bf Step 1.} 
We first settle the continuity problem along $\Lambda_{\rm degen}$. 

\begin{prop}\label{prop:degen}
If $\lambda\in \Lambda_{\rm degen}$ is close enough to zero, then
$f_{\lambda}$ reduces to a 1-dimensional map with 
 entropy $\log d$, and $\chi^+(f_\lambda)\cv
\chi^+(f_0)= \chi^+(p)$ as $\lambda\cv 0$ along $\Lambda_{\rm degen}$
\end{prop}

\begin{proof} We use formulas \eqref{eq:bs5} and 
\eqref{eq:manningbis} --it is possible to
  give a direct proof in the spirit of \cite{mane},  leading 
  however to
  delicate potential theoretic estimates. 

Let $\lambda\in \Lambda_{\rm degen}$ be close to 0, and $L$ be a line
close to $\set{w=0}$. Then by Lemma
\ref{lem:gather2}, $f_\lambda(L)= 
f_\lambda(\cd)$ is a subvariety of some degree $\deg(M_\lambda)$, which is close to
$\set{w=0}$ on compact sets. Let $k=\deg(M_\lambda)$. 
Since $f_\lambda(z,w)$ is close
to $(p(z), 0)$, $M_\lambda$ is a union of graphs over $\set{w=0}$ for large
$\abs{z}$. More precisely if $R>\max\set{p(c),\ c \text{ critical}}$,
then for any simply connected open set $Q\subset\set{\abs{z}>R}$ and
 $\lambda$ small enough (depending on $Q$)
$M_\lambda\cap \pi_1^{-1}(Q)$ is the union of
 $k$ graphs over $Q$
 --see below Lemma
\ref{lem:tang} for more details on a similar argument. 

\medskip

For $\lambda\in \Lambda_{\rm degen}$, we
can define a critical measure associated to the current
$T^-_\lambda$ and the function $G^+_\lambda$, 
exactly in the same way as in Definition \ref{def:critmeas}. Of course
for $\lambda=0$ we  obtain $\mu_{c,0}$. For general $\lambda$,
however, there is an extra $1/k$ factor against each
Dirac mass, coming from the normalization of $T^-_\lambda$. 

In virtue of the description of $M_\lambda$ given above and Lemma
\ref{lem:gather}, we infer that if $A$ is large enough and chosen
so that $\set{G_p=A}$ avoids the postcritical set of $p$, 
$$\int_{A\leq G^+_\lambda<dA} G^+_\lambda d\mu_{c,\lambda} \longrightarrow 
\int_{A\leq G_p<dA} G_p d\mu_{c,0} \text{, as } \lambda\cv 0
\text{ along } \Lambda_{\rm degen}.$$
 
The only  remaining issue is to show that $\log d + 
\int_{A\leq G_\lambda<dA} G^-_\lambda d\mu_{c,\lambda}$ is indeed the
Lyapunov exponent of $f_\lambda$.  

For this, let $\phi=(\phi_1, \phi_2): \cc\cv M_\lambda$ be a
parameterization. Outside the singular set of $M_\lambda$, $\phi$ is
1-1, so $\phi_1:\cc\cv\cc$ is a holomorphic map of degree
$k$, that is, a polynomial of degree
$k$. Normalize so that $\phi_1(t)=t^{k} +
l.o.t.$; the global map $f_\lambda$ induces a polynomial map
$p_\lambda$ on the $t$
variable, in the sense that  
$f_\lambda\rond \phi = \phi\rond p_\lambda$. 
With the normalization done, $p_\lambda$ is monic. Let $G_\phi$
be its Green function in the $t$ coordinate. Since $G^+_\lambda(\phi(t))
= k\log\abs{t} + O(1)$ at infinity, and $G^+_\lambda\rond\phi$ satisfies
the same functional equation as $G_\phi$, we get that  $G^+_\lambda\rond\phi=
kG_\phi$. 

In the $t$ coordinate, we can apply formula \eqref{eq:manningbis} in
the fundamental domain $\set{A\leq kG_\phi<dA}$,
which yields 
$$\chi^+(p_\lambda)=  \chi^+(f_\lambda\rest{M_\lambda}) = \log d +
\int_{A\leq kG_\phi<dA}
G_\phi d\mu_{c, \phi} = \log d + \int_{A\leq G^+_\lambda\circ\phi<dA} 
G^+_\lambda \rond\phi \unsur{k}d\mu_{c, \phi} 
 $$ where
$\mu_{c,\phi}$ is the critical measure associated to $p_\lambda$.  
This finishes the proof, because $\phi_*\left(\unsur{k}\mu_{c,
    \phi}\right)$ equals $\mu_{c,\lambda}$. 
\end{proof} 

\medskip

\noindent {\bf Step 2.} 
From now on, for ease of
reading, it is understood
that for $\lambda\neq 0$, $f_\lambda$ is an automorphism, while $f_0$
is the degenerate map, that is, we restrict ourselves to $\Lambda\setminus
\Lambda_{\rm degen}$. 

In the next lemma, which is purely 1-dimensional dynamics, 
 we give a precise counting of ramified inverse
branches outside the filled Julia set $K_p$. 

\begin{lem}\label{lem:branch}
Let $Q$ be an open topological disk lying outside $K_p$, 
 intersecting  every orbit at most
once. Let $n\geq 1$ and 
$c_1, \ldots, c_q$ be the 
critical points falling in $Q$ after at least one and 
at most  $n$ iterations, and
$j_k$ be the unique integer such that $p^{j_k}(c_k)\in Q$. 

Then the sum of the multiplicities of the critical points of 
$p^n:p^{-n}(Q)\cv Q$ is  
\begin{equation}\label{eq:mult}
\sum_{k=1}^q \mult_z(c_k) d^{n-j_k}.
\end{equation}
\end{lem}

\begin{proof}
The preimage $p^{-n}(Q)$ is a union of topological
  disks. The Riemann Hurwitz formula asserts that 
 the number of these disks is precisely $d^n$ minus 
  the sum of the multiplicities of the critical points of 
$p^n:p^{-n}(Q)\cv Q$.

 Recall the notation $Q^\delta= \set{z\in Q,\
 \mathrm{dist} (z,\fr Q)>\delta}$. Since $Q$ lies outside $K_p$, there
 exists $\delta>0$ such that $Q\setminus Q^\delta$ does not intersect
 the $n^{\rm th}$ image of the critical set. 
In particular, replacing $Q$ by $Q^\delta$, we can assume that there
are no postcritical points on $\fr Q$. 

The topology of $p^{-n}(Q)$ is thus stable under small perturbations of $p$
in the space of polynomials of degree $d$. In particular, without
affecting the sum in \eqref{eq:mult}, we may always
assume that the $c_k$ are critical points of multiplicity 1 and 
that there are no critical orbit relations between the $c_k$ (i.e. the
grand orbits of the $c_k$ are disjoint).  

\medskip

Now subdivide $Q$ into finitely many topological disks $Q_k$, each of
which intersecting the postcritical set in a single point, and so that
no postcritical point is at the boundary of some $Q_k$. Counting the
sum of multiplicities of  the critical points of 
$p^n:p^{-n}(Q)\cv Q$ boils down to doing the computation for each
$Q_k$ and summing over $k$. 

For a single $Q_k$, \eqref{eq:mult} is obvious. Indeed, for $j\leq j_k$,
$p^{-j}(Q_k)$ consists of $d^j$ univalent preimages of $Q_k$. Exactly one of
the components of $p^{-j_k}(Q_k)$ contains the critical point $c_k$
(which is simple by assumption). So among the components of 
$p^{-n}(Q_k)$, exactly $d^{n-j_k}$ are double branched covers of $Q_k$
under $p^n$, while $p^n$ is univalent on the remaining ones. This
concludes the proof.
\end{proof}

We have the following corollary.

\begin{cor}\label{cor:branch}
Under the assumptions of the previous lemma, assume further that 
$Q\subset \set{A<G_p<dA}$, where $A>G_{\rm max}= \max \set{G_p(c), \
c \text{ critical}}$, and let 
 $N$ be so large that $d^{-N}A<\min\set{G_p(c), \ c \text{
     escaping}}$. 
Then the sum of the multiplicities of the critical points of 
$p^N:p^{-N}(Q)\cv Q$ is $d^N \mu_{c, 0}(Q)$.
\end{cor}

\medskip

{\bf Step 3.} We  turn to 2-dimensional dynamics. 
We fix a large $R$ as in the statement of Lemma \ref{lem:gather}, and such
that moreover $\varphi_\lambda\rest{\set{w=0}}$ is univalent on
$\set{\abs{\zeta}>R-C}$ for every $\lambda\in \Lambda$.
This is possible because of the uniformity of the $O(1)$ in 
$\varphi_\lambda(\zeta, 0) = \zeta +O(1)$. Here $C$ is a constant such
that $\abs{\varphi_\lambda(z,w)-z}<C$ on $V_R^+$ for all $\lambda$.
   
Let $A$ be 
such that $\sup_{D^2_R} G_\lambda < A$ for every $\lambda$, and 
such that  there is no postcritical point of $p$ on
 $\set{G_p=A}$.
 Then fix an integer $N$ such that 
\begin{itemize}
\item[-] $A/d^{N-1}< \min\set{G_p(c),
   c\text{ escaping}}$,
\item[-]  $d^{N-1} G_\lambda > A$ on $\overline{V_R^+}$ for all $\lambda$.
\end{itemize} 
These quantities will be kept fixed throughout the proof. 

As before, we project along the fibers of the B{\"o}ttcher
function. As opposite to the previous section, we do not  use
$\varphi_\lambda$ itself but rather a
projection which is $\pi_1$ when $\lambda=0$. Recall that 
$\varphi_\lambda\rest{\set{w=0}}$ is univalent on $\set{\abs{\zeta}>R}$.
The projection along the fibers of $\varphi_\lambda$ onto $\set{w=0}$ is
$\pi_{1,\lambda}:= (\varphi_\lambda \rest{\set{w=0}})^{-1} \rond
\varphi_\lambda$, which is well defined in $V_R^+$ and converges to $\pi_1$
when $\lambda\cv 0$.

The next lemma  establishes the basic 
 connection between 1- and 2-dimensional critical points.

\begin{lem} \label{lem:tang}
Let $Q$ be an open topological disk in $\set{A<G_p<dA}$, intersecting 
every orbit of $p$ at most once. Let $\delta>0$ be such that
$Q\setminus Q^{3\delta}$ does not intersect the postcritical set of
$p$. Let $N$ be the above fixed integer.

Then there exists a neighborhood $\Lambda'$ of $0\in\Lambda$ such
that if $\lambda \in \Lambda'$ and $L$ is any horizontal line of the form 
$L=\set{w=w_0}$, with $\abs{w_0}\leq R$, then:
\begin{enumerate} 
\item[i.] The total number of tangencies
between $f_\lambda^N(L)$ and the fibers of the form $\pi_{1,
  \lambda}^{-1}(\zeta)$, with 
 $\zeta\in Q^\delta$ (resp. of vertical tangencies of $f_\lambda^N(L)$
over $Q^\delta$), 
counted with multiplicity is exactly $d^N \mu_{c,0}(Q)$. 
\item[ii.] There are no tangencies between $f_\lambda^N(L)$ and the 
fibers $\pi_{1,
  \lambda}^{-1}(\zeta)$ for $\zeta \in Q^\delta \setminus
\overline{Q^{2\delta}}$. 
\end{enumerate}
\end{lem}

\begin{proof} For simplicity 
let us consider genuine vertical tangencies first. From
  the expression \eqref{eq:degen} of $f_\lambda$, we get that 
$$f_\lambda^N(z,w)= (p^N(z)+ P_{\lambda, N}(z,w),  Q_{\lambda, N}(z,w)),$$
where $ P_{\lambda, N}$ and $ Q_{\lambda, N}$ vanish for
$\lambda=0$. If we fix a horizontal line $L=\set{w=w_0}$, with $\abs{w_0}\leq
R$,  then $\pi_1\rond f_\lambda^n\rest{L}$ is the map $z\mapsto p^N(z)+
P_{\lambda}^N(z,w_0)$, which is a close perturbation of $p^N$ if
$\lambda$ is small. 

There are no critical points of $p^N$ in   $p^{-N}(Q\setminus
Q^{3\delta})$ 
so, as in Lemma \ref{lem:branch} above, the total number of
 critical points over $Q^\delta$,
counted with multiplicity, is stable under small perturbations. 
From Corollary \ref{cor:branch} 
we conclude that the number of vertical tangencies of $f_\lambda^N(L)$
over $Q^\delta$ equals $d^N \mu_{c,0}(Q)$ if $\lambda$ is small
enough (assertion {\em i.}). 
 The same argument says that there will be no vertical
tangency over $Q^\delta \setminus \overline{Q^{2\delta}}$
(assertion {\em ii.}).

\medskip

The projection $\pi_{1,\lambda}$ is a small perturbation of $\pi_1$
when $\lambda$ is small. Hence by the same
reasoning as above we conclude that the number of tangencies between
$f_\lambda^N(L)$ and the fibers $\pi_{1,\lambda}^{-1}(\zeta)$, $\zeta\in
Q^\delta$ equals $d^N \mu_{c,0}(Q)$ if $\lambda$ is small. Also there
will be no tangencies over $Q^\delta \setminus
\overline{Q^{2\delta}}$
\end{proof}

After horizontal lines, we now
investigate how $f_\lambda^N(D_R^2\cup \overline{V_R^+})$ wraps around over $Q$,
where $Q$ is as in the previous lemma. 
Observe that by the second condition in the choice of $N$, 
$f_\lambda^N (\overline{V_R^+})$ does not intersect $\pi_{1,\lambda}^{-1}(Q)$, so
considering $f_\lambda^N(D_R^2)$ is enough.  

Let $M$ be a horizontal submanifold relative to some fibration
$\pi:\pi^{-1}(D)\cv D$ (cf. \S \ref{subs:basics}),  
and $U(M)$ be some neighborhood of $M$. We say that 
$U(M)$ is a {\em fiberwise  
trivial neighborhood of $M$  over $D'\subset D$} 
if the  projection $\pi$ makes $U(M)$ a
trivial fibration over $M$ in $D'\times\cc$. In other words, 
$U(M)$ is fiberwise  
trivial over $D'$ if for
$\zeta\in D'$, the intersection of $U(M)$ with the fiber
$\pi^{-1}(\zeta)$ consists of  $\mathrm{deg}(M)$ topological
disks, each of which 
containing a single point of $\pi^{-1}(\zeta)\cap M$. 

\begin{lem}\label{lem:tub} 
Let $Q$, $R$, $A$, $\delta$ , $N$  as in  Lemma \ref{lem:tang}, and
$\lambda\in \Lambda$ be so small 
 that the conclusions of Lemma \ref{lem:tang} hold.
Then if $L\subset D_R^2 $ is any horizontal line, 
$f^N_\lambda(D_R^2 )$ is a neighborhood of 
 $f^N_\lambda(L)$, which is 
fiberwise trivial  over $Q^\delta \setminus \overline{Q^{2\delta}}$,
 relative to the projection $\pi_{1,\lambda}$.
\end{lem}

\begin{proof} Recall from the previous lemma that for every horizontal
  line $L_{w_0}=\set{w=w_0}$ 
in $D_R^2 $, $f^N_\lambda(L)$ has no tangencies with
  the fibers $\pi_{1,\lambda}^{-1}(\zeta)$ for
  $\zeta\in  Q^\delta \setminus
\overline{Q^{2\delta}}$. Fix such a fiber
$F_{\zeta_0}= \pi_{1,\lambda}^{-1}(\zeta_0)$  and let us analyze 
$f^N_\lambda(L_{w_0})\cap F_{\zeta_0}$ when $w_0$ ranges across  $D_R$. 
The intersection consists of $d^N$ points, moving holomorphically and
without collision, because $f_\lambda$ is a diffeomorphism and all
intersections are always transverse. We conclude that 
$f^N_\lambda(D_R^2 )\cap F_{\zeta_0}$ is the union of $d^N$ topological
disks, each of which containing a single point of
$f^N_\lambda(L_{w_0})$, for every $w_0$, which was the assertion 
 to be proved. 
\end{proof}

The following easy lemma will be useful for bounding the topology of the
leaves of $T^-$ from below. 

\begin{lem}\label{lem:holonomy}
Let $M$  be a horizontal submanifold,
relative to some fibration
$\pi:\pi^{-1}(D)\cv D$, with no  tangencies with the fibers of $\pi$
over some neighborhood of $\fr D$. 
Assume that both $M$ and the
basis $D$ are isomorphic to disks. 
Let $U(M)$ be some neighborhood of $M$, fiberwise trivial over a
neighborhood of $\fr D$. 

Let $M'\subset U(M)$ be any submanifold. Then $\mathrm{deg}(M)$
 divides $\mathrm{deg}(M')$. 
\end{lem}

\begin{proof} By assumption, $\pi:M\cv D$ is a covering over an annulus 
of the form $D\setminus D'$. Fix a point $\zeta_0\in D\setminus D'$
and a loop $\gamma$ at $\zeta_0$, generating the fundamental group
of $D\setminus D'$. 

Lifting $\gamma$ to $M$ yields a holonomy map of the fiber $M\cap
\pi^{-1} (\zeta_0)$, which is a cycle on the $\mathrm{deg}(M)$
points of the fiber because $M$ has a single
boundary component.

Now, over some neighborhood of $\fr D$, $\pi$
makes $U(M)$ a locally trivial fibration over $M$, and $M'$ is
horizontal with respect to this fibration. In particular, the number
of points of $M'$ in each component of  $U(M)\cap \pi^{-1}
(\zeta)$, counting multiplicities, is a constant $k$. We conclude that
$M'$ has degree $k\mathrm{deg}(V)$ over $D$. 
\end{proof} 
\medskip

{\bf Step 4: conclusion.} Let $A$ and 
 $N$ be as defined before Lemma
\ref{lem:tang}. There are finitely many postcritical points of $p$ in
$\set{A<G_p<dA}$. For each of these postcritical points, we fix a
topological disk $Q$ containing it, and satisfying the requirements of 
Lemma \ref{lem:tang}. 
We will prove that for every small enough $\lambda$,
$$\liminf_{\lambda\cv 0} \mu_{c,\lambda}
\left(\pi_{1,\lambda}^{-1}(Q^\delta)\right)  \geq \mu_{c,0}(Q).$$
As in the proof of Theorem \ref{thm:main}, this implies the lower
semicontinuity of the exponent. 

Consider one of the disks $Q$. By Lemma \ref{lem:tang}, 
for small $\lambda$, $f^N_\lambda(L)\cap \pi_{1,\lambda}^{-1}(Q)$ 
is the union of finitely many horizontal submanifolds $V_1,\ldots ,V_q$, with
the property that: 
\begin{itemize}
\item[-] the $V_i$ have no tangencies with the fibers of
  $\pi_{1,\lambda}$ over $Q^\delta\setminus \overline{Q^{2\delta}}$;
\item[-] the total number of tangencies over $Q^\delta$ equals
  $d^N\mu_{c,0}(Q)$. 
\end{itemize}
By the Maximum Principle, the $V_i$ are topological disks, so by the
Riemann-Hurwitz formula, we get that 
\begin{equation}\label{eq:RH}
\# \text{tangencies over }Q^\delta= 
d^N\mu_{c,0}(Q) = \sum_{i=1}^q (\mathrm{\deg}(V_i)-1).
\end{equation}

Consider now the current $T^-_\lambda$ in $\pi_{1,\lambda}^{-1}(Q)$,
which is horizontal and of slice mass 1 relative to the projection
$\pi_{1,\lambda}$.  The support of $T^-_\lambda$ is contained in
$D_R^2\cup\overline{V^+_R}$, hence by invariance it is contained in
$f^N_\lambda(D_R^2\cup\overline{V^+_R})$.  Also,
$f^N_\lambda(D_R^2\cup\overline{V^+_R} )\cap \pi_{1,\lambda}^{-1}
(Q)=f^N_\lambda(D_R^2 )\cap \pi_{1,\lambda}^{-1} (Q)$ is a
neighborhood of $f^N_\lambda(L)$, which is fiberwise trivial over
$Q^\delta\setminus \overline{Q^{2\delta}}$ by Lemma \ref{lem:tub}. In
particular, if $\zeta_0\in Q^\delta\setminus \overline{Q^{2\delta}}$,
the intersection of this open set with the fiber
$\pi_{1,\lambda}^{-1}(\zeta_0)$ is the union of $d^N$ topological
disks. We claim that the slice mass of $T_\lambda^-$ in each of these
disks is $d^{-N}$.

Indeed, since the sum is 1, it is enough to bound each slice mass from
below by $d^{-N}$. Let $\Delta$ be any of the 
components of $f^N_\lambda(D_R^2 )\cap 
\pi_{1,\lambda}^{-1} (\zeta_0)$. Then $f_\lambda^{-N}(\Delta)$ is a vertical 
 disk in $D_R^2$, so the mass of its intersection with $T_\lambda^-$ is a non zero
 integer.  By the invariance
 property of $T_\lambda^-$, we thus get that the mass of $T_\lambda^-\wedge [\Delta]$
 is $k d^{-N}$, with $k\geq 1$, which yields the desired result.

\medskip 

As in the proof of Theorem \ref{thm:main}, we now use the
decomposition of $T^-_\lambda$ given by Theorem \ref{thm:structure}.
We need to prove that $f_\lambda$ is unstably disconnected for small
$\lambda$. Indeed $Q$ contains a postcritical point, so
$\mu_{c,0}(Q)>0$. Hence by \eqref{eq:RH}, at least one of the $V_i$
has degree $>1$. Let $U(V_i)$ be the connected component of
$f^N_\lambda(D_R^2 )\cap \pi_{1,\lambda}^{-1} (Q^\delta)$ containing
it. The restriction of $T_\lambda^-$ to $U(V_i)$ has slice mass
$\mathrm{deg}(V_i)d^{-N}$ by the previous claim, so it is non zero. In
particular $K^-_\lambda\cap U(V_i)\neq \emptyset$.  If $f_\lambda$ was
unstably connected, then by \cite{bs6}, $K^-_\lambda$ would be the
support of a lamination by graphs over  $\pi_{1, \lambda}$. 
This clearly contradicts Lemma \ref{lem:holonomy}.

\medskip 

We now give a quantitative version of this argument. Consider any of the
horizontal disks $V_i$, and let $U(V_i)$  as above.
Recall that the restriction of $T_\lambda^-$ to $U(V_i)$ has slice mass
$\mathrm{deg}(V_i)d^{-N}$. By Theorem
\ref{thm:structure}, we have a decomposition relative to the
projection $\pi_{1,\lambda}$ over $Q^\delta$,
$$T_\lambda^-\rest{U(V_i)}= \sum_{k=1}^\infty T_k, \text{ with } 
\sum_{k=1}^\infty \ms(T_k)=\mathrm{deg}(V_i)d^{-N}.$$ Now by Lemma
\ref{lem:holonomy} above, the smallest possible degree in the
decomposition is $\mathrm{deg}(V_i)$. Also the
function  $k\mapsto \frac{k-1}{k}$ is increasing.
Thus by Proposition \ref{prop:mass}, we can
bound the critical mass of $T^-_\lambda\rest{U(V_i)}$ from below: 
$$\mu_{c,\lambda}(U(V_i))= \sum_{k=\deg(V_i)}^\infty \frac{k-1}{k}\ms(T_k)
\geq \frac{ \mathrm{deg}(V_i)-1}{
  \mathrm{deg}(V_i)} \sum_{k=\deg(V_i)}^\infty \ms(T_k) = 
 (\mathrm{deg}(V_i)-1) d^{-N}.$$
By summing over $i$, we conclude that the critical mass of
$\pi_{1,\lambda}^{-1}(Q^\delta)$ is bounded from below by 
$\sum_i(\mathrm{deg}(V_i)-1) d^{-N} = \mu_{c,0}(Q)$.
This concludes the proof. 
\end{proof}


\appendix
\section{The fastest escaping critical point}

For the purpose of studying parameter families of polynomial
automorphisms of $\cd$, it is useful to have natural dynamically
defined functions on parameter space. Lyapunov exponents of the
maximal entropy measures are such functions. Here we define a 
notion of ``fastest  rate of escape for critical points'', which is a natural
generalization of 
$G_{\rm max}=\max\set{G_p(c),\ c \text{ escaping}}$ when $p$ is
a polynomial in $\cc$.  In the space of polynomials, this defines a 
psh function which  plays an important role
in \cite{df}. Observe also that $\chi^+(p)\leq \log d + (d-1)G_{\rm
  max}$ as it easily follows from the Manning-Przytycki formula.

So let $f$ be a polynomial automorphism of $\cd$, normalized so that
$f(z,w))= (z^d, 0)+ l.o.t.$ Recall that if $U$ is any open set
avoiding $K^+$, due to the functional equation for $\varphi^+$, 
there exists an integer $N$ so that $(\varphi^+)^N$ is
well defined on $U$. Our definition of $G^+_{\rm max}(f)$ is inspired by
1-dimensional dynamics. 

\begin{defi}\label{def:fastest}
We define the fastest escape rate 
 $G^+_{\rm max}(f)$ as the infimum of $R>0$ such that 
there exists an extension of $\varphi^+$ to a 
neighborhood of $K^-\cap \set{G^->R}$.
\end{defi}

For instance, Bedford and Smillie prove in \cite{bs6} that if $G^+_{\rm
  max}= 0$, then $f$ is unstably connected, hence has no unstable
  critical points. It is likely that $f\mapsto G^+_{\rm max}(f)$ is psh, but
  we could not prove it.
 
As in dimension 1, there is an upper estimate for $\chi^+$ using
$G^+_{\rm max}$. This  should be compared with the estimate given for
horseshoes in \cite[Theorem A.4]{bs5}\footnote{We believe that there
  is a $1/d$ missing in the estimate of $\Lambda-\log d$ given there.}.

\begin{thm}\label{thm:appendix}
$\displaystyle{\chi^+(f)\leq \log d+ dG^+_{\rm max}(f)}$
\end{thm}

\begin{proof} The estimate is based on Proposition \ref{prop:mass}. 
Let $R>G^+_{\rm max}$ be such that the critical measure puts  no mass on
$\set{G^+=R}$. Then $\varphi^+$ is well defined on some neighborhood of
$K^-\cap \set{R< G^+<dR}$, and maps onto the annulus 
$A=\set{e^R<\abs{\zeta}<e^{dR}}$. As before we
 write $\varphi$ for $\varphi^+$.
Then $$\chi^+(f)= \log d + \int_{R\leq G^+<dR} G^+ d\mu_c =
\log d + \int_{\varphi^{-1}(A)} G^+d\mu_c. $$ We claim that 
$\mu_c^-\left(\varphi^{-1}(A)\right)\leq 1$, which
 implies that $\chi^+(f)\leq \log d + dR$. Since $R$ can be arbitrary
 close to $  G^+_{\rm max}$, this will finish the proof.

 \medskip
 
 The first claim is that in $\varphi^{-1}(A)$, $T^-$ admits a
 decomposition $T^-= \sum T_k$ where $T_k$ is made of submanifolds of
 degree $k$ relative to $\varphi$. Indeed, by the discussion preceding
 Proposition \ref{prop:mass}, we know that such a decomposition exists
 in $M\leq G^+ <dM$ for large $M$. Since the existence of such a
 decomposition is invariant under the diffeomorphism $f$, we get the
 result in the original annulus by iterating sufficiently many times. 

Now consider a  neighborhood of $K^-\cap \set{R<G^+<dR}$ of the form 
 $\set{R<G^+<dR}\cap\set{G^-<\e}$ where $\varphi$ is well defined. The
 projection $\varphi$ needn't be a locally trivial fibration 
onto $A$, however we leave the reader check that the assumption of
 Remark \ref{rmk:appendix} is satisfied.  In particular, we can say 
that $T^-$ and the 
 $T_k$ are horizontal currents relative to the projection $\varphi$
 and the slice mass is invariant. 

Cut the annulus $A$ by a radial line $L$ such that $\varphi^{-1}(L)$
has zero  mass for the critical measure. Then $Q:=A\setminus L$ is a
simply connected open set, and we have a decomposition (different from
the previous one) $T^-=\sum T_k$ of $T^-$ over $Q$ relative to
$\varphi$. By  Proposition \ref{prop:mass}, we infer that 
$$\mu_c(\varphi^{-1}(Q)) = 
\sum_{k=1}^\infty \frac{k-1}{k} \ms(T_k).$$
Near infinity, the slice mass of $T^-$ relative to $\varphi$ is 1, so
by invariance of the slice mass, this is also true in $A$. We conclude
that $\sum_{k=1}^\infty  \ms(T_k)=1$,
 hence  $\mu_c(\varphi^{-1}(Q))\leq 1$. 
\end{proof}


\begin{thebibliography}{[ABCD]}

\bibitem[BB]{bas-ber} 
Bassanelli, Giovanni; Berteloot, Fran{\c c}ois.
\newblock \emph{Bifurcation currents in holomorphic dynamics on ${\bf
    P}^k$.}
\newblock   Preprint (2005), {\tt math.DS/0507555}. To
appear in J. Reine Angew. Math.

\bibitem[BJ]{bj} Bedford, Eric; Jonsson, Mattias.
\newblock {\em Dynamics of
    regular polynomial endomorphisms of $\cc^k$}.
\newblock Amer. J. Math.
122 (2000), 153-212.

\bibitem[BS1]{bs1} Bedford, Eric; Smillie, John. 
\newblock {\em Polynomial
diffeomorphisms of $\cc^ 2$: currents, equilibrium measure and
hyperbolicity.} 
\newblock Invent. Math. 103 (1991), 69-99. 

\bibitem[BS3]{bs3} Bedford, Eric; Smillie, John. 
\newblock {\em Polynomial
diffeomorphisms of $\cc^ 2$. III. Ergodicity, exponents and entropy of
the equilibrium measure.} 
\newblock Math. Ann. 294 (1992),  395-420.

\bibitem[BLS1]{bls}
Bedford, Eric;  Lyubich, Mikhail;  Smillie, John.
\newblock  {\em Polynomial
diffeomorphisms of $\cc^ 2$. IV: The measure of maximal entropy and
laminar currents. }
\newblock Invent. Math.  112  (1993), 77-125.

\bibitem[BLS2]{bls2}
Bedford, Eric;  Lyubich, Mikhail;  Smillie, John.
\newblock{\em Distribution of periodic points of polynomial diffeomorphisms of
  $\cd$.}
\newblock Invent. Math.  114 (1993), 277-288.

\bibitem[BS5]{bs5} Bedford, Eric; Smillie, John. 
\newblock{\em Polynomial diffeomorphisms of $\cd$. V. Critical points and
  Lyapunov exponents.} 
\newblock  J. Geom. Anal.  8  (1998),  349-383.

\bibitem[BS6]{bs6} Bedford, Eric; Smillie, John.
\newblock{\em Polynomial diffeomorphisms of $\cd$. VI. Connectivity of  $J$.}
\newblock  Ann. of Math. (2)  148  (1998),   695-735. 

\bibitem[Bi]{bishop} Bishop, Errett
\newblock{\em Conditions for the analyticity of certain sets.}
\newblock  Michigan Math. J.  11  (1964) 289-304.

\bibitem[Dem]{de}
Demailly, J.-P.
\newblock\emph{Monge-Amp{\`e}re operators, Lelong numbers and intersection 
theory.}
\newblock Complex analysis and geometry,  115--193, Univ. Ser. Math., 
Plenum, New York, 1993. 


\bibitem[DeM]{demarco2}
DeMarco, Laura.
\newblock\emph{Dynamics of rational maps: Lyapunov exponents, 
bifurcations, and capacity.}
\newblock  Math. Ann.  326  (2003), 43-73.

\bibitem[DS]{ds}
Dinh, Tien Cuong; Sibony, Nessim.
\newblock{\em Geometry of currents, intersection theory and dynamics of
  horizontal-like maps.} 
\newblock  Ann. Inst. Fourier (Grenoble)  56 (2006),  423-457.

\bibitem[Du1]{hl} Dujardin, Romain. 
\newblock{\em H{\'e}non-like mappings in    $\cd$.} 
\newblock Amer. J. Math. 126 (2004), 439-472. 

\bibitem[Du2]{structure} Dujardin, Romain. 
\newblock{\em Structure properties
    of laminar currents on $\pd$}. 
\newblock  J. Geom. Anal.  15  (2005),  25-47.

\bibitem[Du3]{connex}  Dujardin, Romain.
\newblock{\em Some remarks on the connectivity of Julia sets for
  2-dimensional diffeomorphisms,}
\newblock in {\em  Complex Dynamics,} 63-84, Contemp. Math., 396, 
Amer. Math. Soc., Providence, RI, 2006.

\bibitem[DF]{df} Dujardin, Romain; Favre, Charles.
\newblock{\em Distribution of rational maps with a preperiodic
  critical point.}
\newblock Preprint (2006), {\tt math.DS/0601612}.

\bibitem[FLM]{flm}
 Freire, Alexandre; Lopes, Artur; Ma{\~n}{\'e}, Ricardo. 
\newblock{\em An invariant measure for rational maps.} 
\newblock Bol. Soc. Brasil. Mat.  14  (1983),   45-62.

\bibitem[FM]{fm} Friedland, Shmuel; Milnor, John.
\newblock{\em Dynamical properties of plane polynomial automorphisms.}
Ergodic Theory Dynam. Systems 9 (1989),  67-99.

\bibitem[HO]{ho} Hubbard, John H.; Oberste-Vorth, Ralph W. 
{\em H{\'e}non mappings in the complex domain. I. The global topology of
  dynamical space.}  Inst. Hautes {\'E}tudes Sci. Publ. Math. No. 79
(1994), 5-46. 

\bibitem[L]{lyub}
Lyubich, Mikhail. 
\newblock \emph{Entropy properties of rational endomorphisms of 
the Riemann sphere.}
\newblock  Ergodic Theory Dynam. Systems  3  (1983),   351-385.

\bibitem[Ma1]{mane}  Ma{\~n}{\'e}, Ricardo.
\newblock{\em The Hausdorff dimension of invariant probabilities of
  rational maps.}
  Dynamical systems, Valparaiso 1986,  86-117, Lecture Notes in Math.,
  1331,  Springer, Berlin, 1988. 

\bibitem[Ma2]{manning} Manning, Anthony. 
\newblock{\em The dimension of the maximal measure for a polynomial
  map. }
\newblock  Ann. of Math. (2)  119  (1984),  425-430.

\bibitem[Ph]{pham}  Pham, Ngoc-mai.
\newblock \emph{Lyapunov exponents and bifurcation current for
  polynomial-like maps.} 
\newblock Preprint (2005), {\tt math.DS/0512557}.

\bibitem[Pr]{prz} Przytycki, Feliks.
\newblock{\em Hausdorff dimension of harmonic 
measure on the boundary of an attractive basin for a holomorphic map.}
\newblock  Invent. Math.  80  (1985),   161-179. 

\bibitem[Si1]{sib}  Sibony, Nessim. 
\newblock{\em Dynamique des applications rationnelles de $\mathbb{P} ^k$}.
\newblock Dynamique et g{\'e}om{\'e}trie complexes (Lyon, 1997), 
  Panoramas et Synth{\`e}ses, 8, 1999. 

\bibitem[Si2]{sib-ucla}  Sibony, Nessim.
\newblock{\em Iteration of polynomials.}
\newblock UCLA Lecture notes (unpublished).


\bibitem[S{\l}]{sl} S{\l}odkowski, Zbigniew.
\newblock {\em
Uniqueness property for positive closed currents in $\cc\sp 2$.}
\newblock Indiana Univ. Math. J. 48 (1999), 635-652.

\bibitem[Y]{young}
Young, Lai Sang 
\newblock {\em Dimension, entropy and Lyapunov exponents.} 
Ergodic Theory Dynamical Systems  2  (1982),  109-124.

\end{thebibliography}
\end{document}